\documentclass [11pt]{amsart}
\usepackage {amsmath, amssymb, amscd, graphicx, color, url,epsf}
\usepackage[all,knot,arc]{xy}
\usepackage[text={6.5in,9in},centering,letterpaper,dvips]{geometry}
\usepackage{hyphenat}

\newcommand\sign{\mathrm{Sign}}
\newcommand\si{\mathrm{sign}}
\newcommand\Short{\mathrm{short}}
\newcommand\Long{\mathrm{long}}

\newcommand\argmax{\arg\max}
\newcommand\rect{\mathrm{Rect}}
\newcommand\bigon{\mathrm{Bigon}}

\newtheorem {theorem}{Theorem}[section]
\newtheorem {lemma}[theorem]{Lemma}

\theoremstyle{definition}
\newtheorem {definition}[theorem]{Definition}

\theoremstyle{remark}

\usepackage{epsf}

\begin{document}

\title{Effective computation of knot Floer homology}

\author{Jean-Marie Droz}
\thanks {JMD is supported by a SNF grant}
\address {Institut f\"ur Mathematik, Universit\"at Z\"urich, Winterthur Strasse 190, CH-8057 Z\"urich Switzerland}
\email {jdroz@math.unizh.ch}

\begin {abstract} 
Beliakova gave in \cite{beliakova} a method for computing link Floer homology. We extend her approach to $\mathbb{Z}$ coefficients and implement it in a publicly available computer program. We review the main programming and optimization methods used. Our program is then used to check that the Floer homology of a prime non-alternating knot with $\leq 11$ crossings has no torsion. 
\end {abstract}

\maketitle

\section{Introduction}
 
Knot Floer homology $HL$ was introduced in \cite{knots} by Ozsv{\'a}th and Szab{\'o} and independently in \cite{rasmussen} by Rasmussen. In this article, we will restrict ourselves to the study of $\widehat{HL}$, a simplified version of $HL$.
The knot Floer homology $\widehat{HL}(K)$ of a knot $K$ is the homology of a bigraded complex $(C,\partial)$ with Maslov grading $m$ and Alexander grading $a$. $$\widehat{HL}(K)=\bigoplus_{m,a \in \mathbb{Z}} HL_{m} (K,a)$$ We will work with either ${\mathbb Z}$ or ${\mathbb Z}/2{\mathbb Z}$ as coefficient ring. Knot Floer homology, even in its simplest version $\widehat{HL}(K)$ over ${\mathbb Z}/2{\mathbb Z}$, gives much information about a knot: The Seifert genus $g(K)$ of $K$ is
$$\underset{a}\argmax (\, \widehat{HL}(K,a)\neq 0)$$ (see \cite{os}).
The knot $K$ is fibered if and only if
$\mathrm{rank\,} \widehat{HL}(K,g(K))=1 $ (see \cite{Ghiggini} and \cite{YiNi}). It is also possible to extract from $HL$ a bound for the slice genus of the knot, see \cite{4BallGenus}.

The first combinatorial way for calculating knot Floer homology was given in \cite{oms} for ${\mathbb Z}/2{\mathbb Z}$ coefficients and was then extended for $\mathbb{Z}$ coefficients in \cite{moz}. Those papers explain how to construct a complex $(C,\partial)$ from the grid diagram (also called rectangular diagram or arc presentation) of $K$, so that $\widehat{HL}(K)=H_{\ast}(C)\otimes  V^{\otimes(n-1)}$ (see Theorem \ref{MOSTh} for the definition of the free module $V$). However, the number of generators in $C$ increases very fast with the size of the grid diagram, making the practical computation of this complex difficult. A computer program extracting knot Floer homology  from this complex was written by Baldwin and Gillam (see \cite{bg}). This program manages to alleviate the problem of the number of generators by taking advantage of the sparseness of the matrix describing the boundary map $\partial$.

In \cite{beliakova}, Beliakova proposed a new complex $(C_{\Long},\partial_{\Long})$ for knot Floer homology, which we will call the long oval complex. She explained how it is homotopic to a complex $(C_{\Short},\partial_{\Short})$, which we will call short oval complex. The main interest of the construction is that the short oval complex usually has a much smaller number of generators than the MOS complex. Those two complexes are described in Section 2. 

In Section 3, we extend Beliakova's complexes, originally defined with ${\mathbb Z}/2{\mathbb Z}$ coefficients, to $\mathbb Z$ coefficients. In Section 4, the formula used to calculate the Alexander grading of generators in oval complexes is proved.

 Section 5 explains how our program computes the short oval complex and extracts its homology. It also describes algorithms and optimizations that could have a larger interest. Many could, for example, be used to work with the MOS complex. 

\subsection{Main results.} We developed a program that can be used to determine the $\widehat{HL}$ with $\mathbb{Z}$ or ${\mathbb Z}/2{\mathbb Z}$ coefficients of knots with less than 13 crossings efficiently and can be used to determine fiberedness and the Seifert genus of even bigger knots\footnote{Our program is accessible on our homepage "http://www.math.unizh.ch/assistenten/jdroz" or via the Knot Atlas "http: /katlas.math.toronto.edu/wiki/Main\_Page"}. Using our program, we show that the $\widehat{HL}$ of prime non-alternating knots with less than 12 crossings contains no torsion and checked the $\widehat{HL}$ computations given in \cite{bg}.

\subsection{Acknowledgements.} I would like, above all, to express my gratitude to Anna Beliakova for the original impetus of the project, numerous interesting discussions, advice and ideas. I also would like to thank Dror Bar-Natan for good advice and his support in making my program a part of the Knot Atlas. The resources of the Knot Atlas\footnote{http: /katlas.math.toronto.edu/wiki/Main\_Page} and of the knot data tables of Alexander Stoimenov\footnote{http://www.kurims.kyoto-u.ac.jp/~stoimeno/ptab/index.html} were very useful to me. Some of the figures in this article are courtesy of Anna Beliakova.

\section{Definitions of the complexes}
\subsection{MOS complex}
\begin{figure}[h]
\mbox{\epsfysize=50mm \epsffile{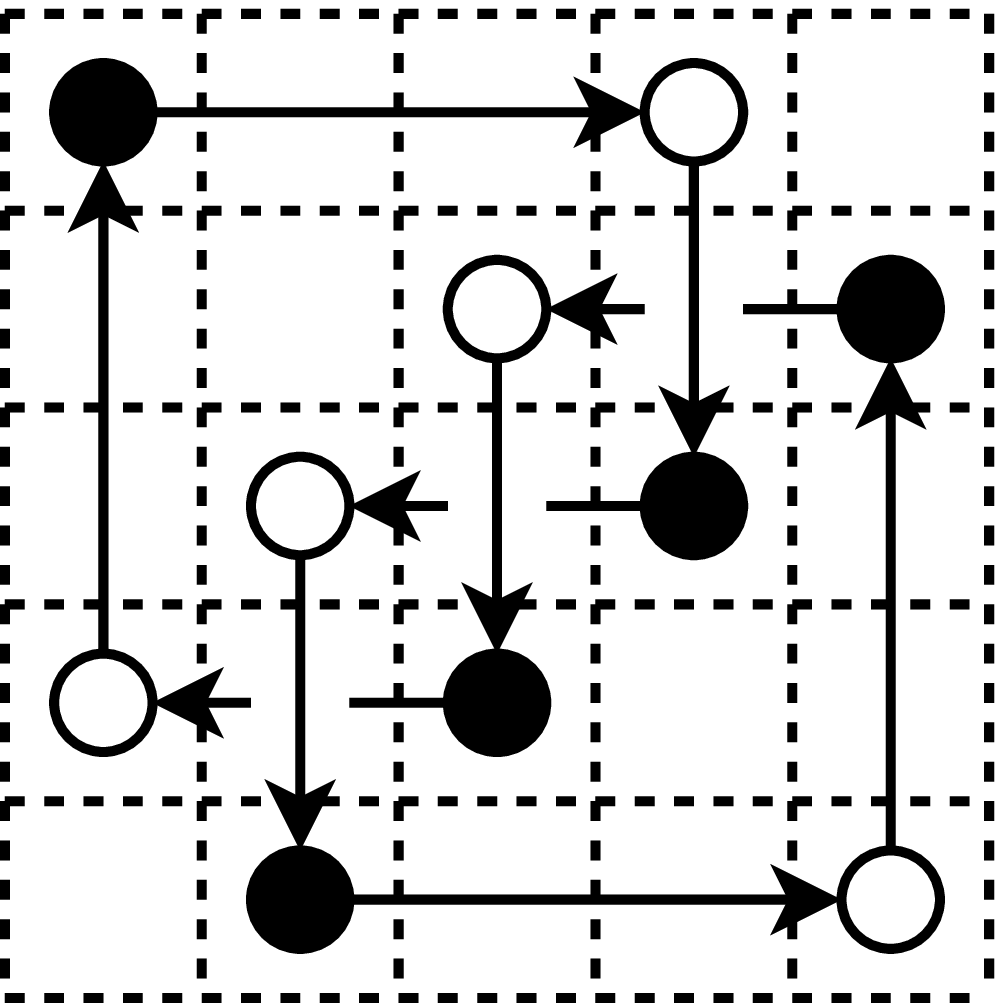}} 

\caption{{\ A grid diagram  for the $3_1$ knot.} The white points represent $\mathbb{X}$ and the black points $\mathbb{O}$.}
\label{gridDiag}
\end{figure}

A \emph{grid diagram} (or grid diagram on the plane) of complexity $n$ is a $n\times n$ square grid on the plane, with each square decorated by an $X$, an $O$ or nothing, such that each column or row of the grid contains exactly one $X$ and one $O$. The set of all $O$ (respectively all $X$) is called $\mathbb{O}$ (respectively $\mathbb{X}$). (See Figure \ref{gridDiag}.)

 We call the center of a decorated square a \emph{puncture}. More information about grid diagrams can be found in \cite{dy}.

A grid diagram represents a link in the following way. If we draw a line segment between pairs constituted of an $X$ and an $O$ sitting in the same row or in the same column, we obtain a figure in the plane. Assuming that vertical segments always pass over horizontal segments when they intersect, we can see this figure as a projection of a link. We say that the grid diagram represents the same link as this projection. It is not difficult to see that a knot can always be represented by many different grid diagrams.

Starting with a grid diagram on the plan, by identifying the uppermost and bottommost lines of the grid and the leftmost and rightmost lines, we obtain a decoration of a toroidal grid called a\emph{ grid diagram on the torus}. After this operation, the vertical and horizontal lines of the grid become meridional and longitudinal circles. Each grid diagram on the torus of complexity $n$ is associated in a natural way to $n^2$ grid diagrams in the plane.

\begin{definition}
Let $D$ be a grid diagram on the torus of complexity $n$, and $D'$ one of its analog in the plane. 
The \emph{MOS complex} $(C(D),\partial(D))$ is defined in the following way: $C(D)$ is the vector space over ${\mathbb{Z}}/2{\mathbb{Z}}$, which has as basis the unordered $n$-tuples of intersection points of circles of the grid that contains exactly one intersection point on every meridional or longitudinal circle. The boundary map is defined as follows: $$\partial(x)=\sum_{y\in \rect_0(x)} y.$$ We define $\rect_0(x)$ to be the set of generators $y$ such that:
\begin{itemize}
\item{The generators $x$ and $y$ have exactly $n-2$ intersections in common.}
\item{The four intersection points that are either only in $x$ or only in $y$ form the four corners of a rectangle $R$. In $D'$, the upper right corner of $R$ is in $x$.}
\item{The rectangle $R$ is empty in the sense that it doesn't contain a puncture of the grid diagram or any intersection point of $x$.}
\end{itemize}
\end{definition} 
We often think of the boundary map as a matrix with entries indexed by the generators of the complex: $$\partial(x)=\sum_{y} \partial_{x,y}\cdot y.$$ In other words, $\partial=(\partial_{x,y})$ for $x,y$ generators.
\begin{definition}
let $S,T \subset \mathbb{R}^2$ be two finite sets of points in the plane. $I(S,T)$ is defined as the number of pairs $(a,b)\in S,\, (c,d) \in T$ with $a<c$ and $b<d$. We also define $J(S,T)=(I(S,T)+I(T,S))/2$. We consider the natural bilinear extensions of $I$ and $J$ to formal sums of intersection points.
\end{definition}

\begin{definition}
\label{masDef}
A generator $g \in C$ has Maslov grading $M(x)=I(x,x)-I(x,\mathbb{O})-I(\mathbb{O},x)+I(\mathbb{O},\mathbb{O})+1$. It will also sometimes be practical to speak of the Maslov grading of an arbitrary set of points.
\end{definition} 
\begin{definition}
\label{alexDef}
A generator $g \in C(D)$ for a grid diagram $D$ of complexity $n$ has Alexander grading $A(x)=J(x-{(\mathbb{O}+\mathbb{X})}/2 ,\mathbb{X}-\mathbb{O})-{(n-1)}/2$. 
\end{definition} 
It is also possible to write the Alexander grading of a generator $x$ in a diagram $D$ as the sum of the winding number of $D$ around the intersection points $p \in x$ and a constant depending only on $D$. Let $a(p)$ denote the average of the winding number of the knot projection represented by $D$ in a small ball around $p$.
$$
\label{AGrading}
A(x)=\sum_{p\in x} a(p)- \frac{1}{2} \Bigl(\sum_{o\in O} 
a(o)\Bigr) - 
\left(\frac{n-1}{2}\right)
$$
 Despite a small difference in the definition of $a(p)$, the formula above is equivalent to the formula for the Alexander grading given in \cite{oms}.

\begin{theorem}
\label{MOSTh}
(C.Manolescu, P.Ozsv\' ath, S.Sarkar) Let $D$ be a grid diagram of complexity $n$ of the knot $K$, $(C(D),\partial(D))$ the MOS complex constructed from $D$.
\begin{itemize}
\item{The chain complex $(C(D),\partial(D))$ is bigraded by the Alexander and Maslov gradings. The boundary map preserves the Alexander grading and decreases the Maslov grading by 1.}
\item{The homology of $H_* (C(D),\partial)$ is $\widehat{HL}(K)\otimes  V^{\otimes(n-1)}$}, where $V$ is a vector space with basis composed of one vector of Alexander and Maslov gradings -1 and one vector with Maslov and Alexander gradings zero.
\end{itemize}
\end{theorem}
\begin{figure}[t]
\mbox{\epsfysize=8cm \epsffile{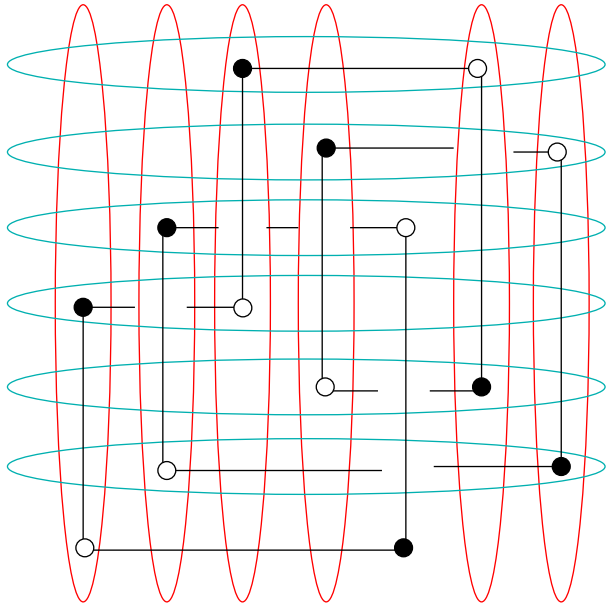}} 
\caption{{\ A grid diagram  for the $5_2$ knot with a collection of long ovals.}}
\label{longFig}
\end{figure}
\subsection{Beliakova's complexes}
\label{BelComp}
We will now construct two complexes that have the same homology as the MOS complex, and for which a theorem almost identical to Theorem \ref{MOSTh} is proven in \cite{beliakova}.
The construction of the \emph{long oval complex} $(C_{\Long}(D),\partial_{\Long}(D))$ also starts with the grid diagram $D(K)$ of the knot $K$. But this time, the construction is done in the plane. Let us call $n$ the complexity of the grid diagram $D$. We begin by drawing $n-1$ long thin vertical (respectively horizontal) ovals around all but one pair of punctures with identical first (respectively second) coordinates (see Figure \ref{longFig}). We call the vertical ovals $\boldsymbol{\alpha}=\{\alpha_1,\ldots ,\alpha_{n-1} \}$ and the horizontal ovals $\boldsymbol{\beta}=\{\beta_1,\ldots ,\beta_{n-1}\}$. While drawing those ovals, we also require that at least one puncture is in the unbounded component of $\mathbb{R}^2/(\alpha_1 \cup ... \alpha_{n-1}\cup \beta_1 \cup ... \beta_{n-1})$.

\begin{definition}
\label{longOvalDef}
Let $D$ be a grid diagram of complexity $n$ in the plane. The vector space over ${\mathbb{Z}}/2{\mathbb{Z}}$ called $C_{\Long}(D)$ has as basis elements the sets of $n-1$ intersection points between vertical and horizontal ovals with exactly one intersection point on each oval.  An Alexander grading on those generators is defined by exactly the same formula as in the case of the MOS complex. The Maslov grading $M(x)$ of a generator $x$ is defined to be $I(x,x)-I(x,\mathbb{O})-I(\mathbb{O},x)+I(\mathbb{O},\mathbb{O})$. (The same formula as for the MOS complex minus one).
 The boundary map $\partial_{\Long}(D)$ is defined as follows: 
$$\partial(x)=\sum_{y\in \rect_0(x)\cup y\in \bigon_0(x)} y.$$
 The set $\rect_0(x)$ contains the generators $y$ such that:
\begin{itemize}
\item{The generators $x$ and $y$ have exactly $n-3$ intersection points in common.}
\item{The four intersection points that are either only in $x$ or only in $y$ form the four corners of a rectangle $R$. The sides of $R$ are arcs of ovals. The upper right corner of $R$ is in $x$.}
\item{The rectangle $R$ is empty in the sense that it doesn't contain any puncture of the grid diagram or any intersection point of $x$.}
\end{itemize}
The set $\bigon_0(x)$ contains the generators $y$ such that:
\begin{itemize}
\item{The generators $x$ and $y$ have exactly $n-2$ intersection points in common.}
\item{The two intersection points that are either only in $x$ or only in $y$ constitutes the two corners of a bigon $L$. The sides of $L$ are one arc of horizontal oval and one arc of vertical oval. A counterclockwise rotation around $L$, along the arc of the horizontal oval, leads from the corner in $x$ to the corner in $y$.}
\item{The bigon $L$ is empty in the sense that it doesn't contain any puncture of the grid diagram.}
\end{itemize}
\end{definition}
\begin{figure}[t]
\mbox{\epsfysize=8cm \epsffile{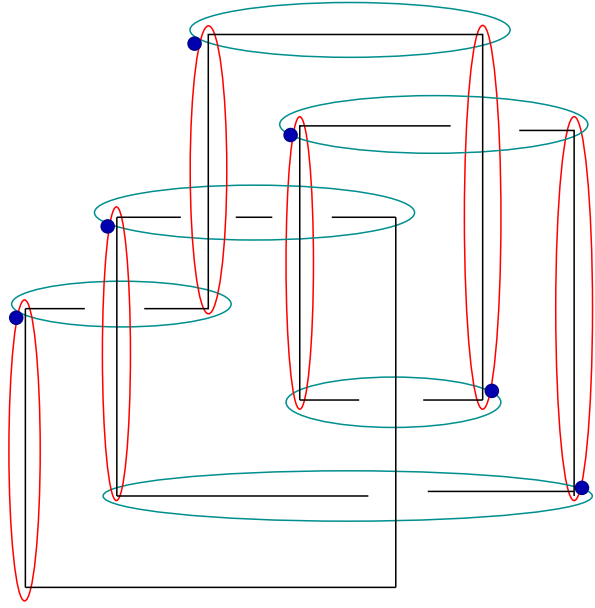}} 

\caption{{\ A grid diagram  for the $5_2$ knot with a collection of short ovals.} The blue points represent a generator.}
\label{shortFig}
\end{figure}

\begin{figure}[h]
\mbox{\epsfysize=40mm \epsffile{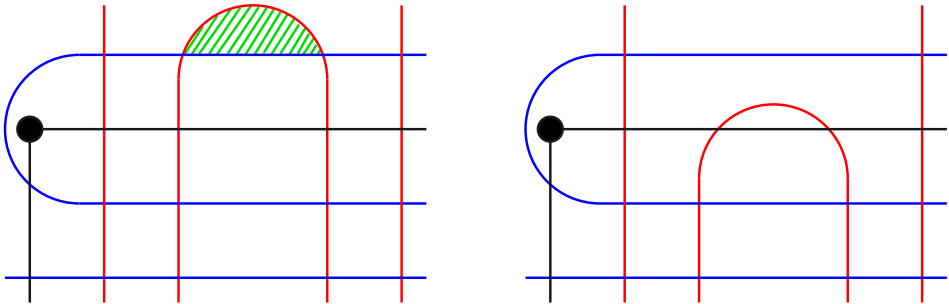}} 
\caption{{\ The disappearance of a bigon during an homotopy from the long oval complex to the short oval complex.}}
\label{bigon2}
\end{figure}

The long oval complex can be reduced, by a sequence of homotopies $(h_1,...,h_l)$, to a much smaller (in terms of the rank of $C_{\Long}(D)$) complex, which we will call the \emph{short oval complex} $(C_{\Short}(D),\partial_{\Short}(D))$ (see Figure \ref{shortFig}). This sequence of homotopies corresponds to a progressive shortening of the curves $\boldsymbol{\alpha}$ and the curves $\boldsymbol{\beta}$. This shortening is an homotopy written $(\boldsymbol{\alpha},\boldsymbol{\beta})_t$ for $t \in [0,l+1]$ such that there is no $t$ for which a curve meets a puncture. Among such homotopies, we choose one that minimizes the number of intersection points of the curves at $t=l+1$. We choose this progressive shortening such that one pair of intersection points disappears for each $t\in \{1,...,l\}$ (see Figure \ref{bigon2} for an example of the corners of a bigon disappearing). To each $(\boldsymbol{\alpha},\boldsymbol{\beta})_t$ $t\in [0,l+1]$, we associate a complex $(C_{\lfloor t\rfloor}(D),\partial_{\lfloor t\rfloor}(D))$. The complex $C_{k}(D)$ is generated by sets of $n-1$ intersection points of $\boldsymbol{\alpha}_k$ and $\boldsymbol{\beta}_k$ such that each intersection point is on exactly one oval. The boundary maps $\partial_k(D)$ are defined inductively. 

\begin{definition}
\label{edgeDel}
Let $p_1$ and $p_2$ be the intersection points that disappear at time $t=k$, assume $\{ p_1\}$ has a bigger Maslov grading than $\{ p_2\}$. 
For $x$ a generator of $C_{k}(D)$, let $\eta (x)$ be either $0$, if $x$ does not contain $p_2$, or, if $x$ does contain $p_2$, a generator identical to $x$ except that it contains $p_1$ instead of $p_2$.
$$\partial_k=\pi \circ(\partial_{k-1}+\partial_{k-1}\circ \eta \circ \partial_{k-1})\circ \iota$$
Where $\pi$ is the natural projection from $C_{k-1}(D)$ to $C_k(D)$ that sends generators containing $p_1$ or $p_2$ to zero and $\iota$ is the natural injection from $C_k(D)$ to $C_{k-1}(D)$. 
\end{definition}

The homotopy $h_k$ is an homotopy between $(C_{k-1}(D),\partial_{k-1}(D))$ and $(C_{k}(D),\partial_{k}(D))$. It corresponds to the "cancellation" of two intersection points. The existence of these homotopies is Lemma 2.1 of \cite{beliakova}.

It is noteworthy that this construction of $(C_{\Short}(D),\partial_{\Short}(D))$ is not canonical. The differential $\partial_{\Short}$ depends on the whole $(\boldsymbol{\alpha},\boldsymbol{\beta})_t$ (for $t\in [0,l+1]$) and not just on $(\boldsymbol{\alpha},\boldsymbol{\beta})_{l+1}$.
\begin{theorem}
\label{BTh}
(Beliakova) Let $D$ be a grid diagram of complexity $n$ of the knot $K$, $(C_{\Short}(D),\partial_{\Short}(D))$ a short oval complex constructed from $D$.
\begin{itemize}
\item{The chain complex $(C_{\Short}(D),\partial_{\Short}(D))$ is bigraded by the Alexander and Maslov gradings. The boundary map preserves the Alexander grading and decreases the Maslov grading by 1.}
\item{The homology $H_* (C_{\Short}(D),\partial_{\Short}(D))$ is $\widehat{HL}(K)\otimes  V^{\otimes(n-1)}$}, where $V$ is a vector space with basis composed of one vector of Alexander and Maslov gradings -1 and one vector with Maslov and Alexander gradings zero.
\end{itemize}
\end{theorem}

\section{Sign Assignment}
We start by defining an extension $(C'_{\Long}(D),\partial '_{\Long}(D))$ of the long oval complex, the homology of which is $\widehat{HL}$ over $\mathbb{Z}$. The complex $(C'_{\Long}(D),\partial '_{\Long}(D))$ will be said to be a sign assignment over the long oval complex. 

\begin{definition}
\label{signAssignment}
Let $(C',\partial')$ be a complex with $\mathbb{Z}$ as coefficient ring and basis $\boldsymbol{b'}=\{b'_1,\ldots,b'_n\}$. Let $(C,\partial)$ be a complex with ${\mathbb{Z}}/2{\mathbb{Z}}$ as coefficient ring and basis $\boldsymbol{b}=\{b_1,\ldots,b_n\}$.
Let us assume that the matrices $(\partial'_{i,j})$ and $(\partial_{i,j})$ represent the boundary maps in the basis $\boldsymbol{b'}$ and $\boldsymbol{b}$. We say that $(C',\partial')$ is a \emph{sign assignment} on $(C,\partial)$ if the following two conditions are true for all integers $1\leq i,j\leq n$.
\begin{itemize}
\item{$\partial_{i,j}=0\: \Rightarrow \: \partial'_{i,j}=0$}
\item{$\partial_{i,j}\neq 0\: \Rightarrow \: \partial'_{i,j}=\pm 1$}
\end{itemize}
\end{definition}

Our approach is completely analogous to the one taken in \cite{moz} for the same purpose. We assume for the following definitions that the curves in $\boldsymbol{\alpha}\cap\boldsymbol{\beta}$ are ordered. We say that a point sits positively on an oval if it sits on the right side of an $\alpha$ curve or the upper side of a $\beta$ curve.
\begin{definition}
\label{signDef}
Let $x$ and $y$ be generators of the long oval complex. We define $\sign(x,y)$ in the following way:
\begin{itemize}
\item{If $y \in \rect_0(x)$ (we say that $x$ and $y$ are \emph{connected by a rectangle}):
Let $(a,b)$ (respectively $(c,d)$) denote the coordinates of the leftmost (respectively rightmost) point that is in $x$ but not in $y$. Let $\mathrm{D}(x,y)$ be the number of points $p=(p_1,p_2)$ in $x$ such that $a\leq p_1 \leq c$ and $p_2 \leq b$. In other words, $\mathrm{D}(x,y)$ is the number of points below the rectangle between $x$ and $y$.
$$\sign(x,y)=(-1)^{I(x,\{ (x_1,x_2)\in x\mid x_2 \leq d\})+\mathrm{D}(x,y)\cdot(I(x,\{ (x_1,x_2)\in x\mid b<x_2\leq d\})+1)}$$
}
\item{If $y \in \bigon_0(x)$ (we say that $x$ and $y$ are \emph{connected by a bigon}):
Let $E$ be the oval on which $x$ and $y$ have intersection points on opposite sides (the intersection points are placed symmetrically with respect to the great axis of $E$). Let $\mathrm{pre}(x,y)$ be the number of ovals that come before $E$ in the ordering and on which an intersection point of $x$ sits positively.
$$\sign(x,y)=(-1)^{I(x,x)+\mathrm{pre}(x,y)}$$}
\item{Otherwise, $\sign(x,y)$ is undefined.}
\end{itemize}
\end{definition}

\begin{definition}
The complex $C'_{\Long}(D)$ is a free module with $\mathbb{Z}$ coefficients on the same generators as $C_{\Long}(D)$.
 The boundary map $\partial '_{\Long}$ is the linear map equal to $\sum_{y\in \rect_0(x) \cup \bigon_0(x)} \sign(x,y) \cdot y$ for each generator $x$ of $C'$. 
\end{definition}
For generators ($x$,$y$) between which there is a rectangle ($y \in \rect_0(x)$), our formula is taken from \cite{moz}.

\begin{lemma}
\label{d2iszero}
The boundary map $\partial '_{\Long}$ verifies the condition ${\partial '}_{\Long}^2=0$ and thus $(C'_{\Long}(D),\partial '_{\Long}(D))$ is a chain complex.
\begin{proof}
Let $x$ and $z$ be two generators of $C'_{\Long}(D)$ with $M(x)-2=M(z)$. The number of generators $y$ such that $y\in \rect_0(x) \cup \bigon_0(x)$ and $z\in \rect_0(y) \cup \bigon_0(y)$ is either 0 or 2. If those two generators exist, we call them $y_1$ and $y_2$.
The lemma is equivalent to the claim that when $y_1$ and $y_2$ exist, $\sign(x,y1)\cdot\sign(x,y2)\cdot\sign(y1,z)\cdot\sign(y2,z)=-1$.
We must check this in three cases:
\begin{itemize}
\item{First case: $y_1,y_2\in \rect_0(x)$. This implies that $z\in \rect_0(y_1)$ and $z\in \rect_0(y_2)$. Therefore, the claim follows from \cite{moz} Section 4.1, where an analogous claim for the MOS complex is proven.}

\item{Second case: $y_1\in \rect_0(x)$ and $y_2 \in \bigon_0(x)$. This implies $z\in \rect_0(y_2)$ and $z \in\bigon_0(y_1)$. Since the formula for $\sign$ for generators connected by a rectangle only depends on which pairs of ovals have common intersection points, we have $\sign(x,y_1)=\sign(y_2,z)$. Since $\mathrm{pre}(x,y_2)=\mathrm{pre}(y_1,z)$ and $I(x,x)\not\equiv I(y_1,y_1) \pmod{2}$, $\sign(x,y_2)=-\sign(y_1,z)$ and the claim follows.}

\item{Third case: $y_1,y_2 \in \bigon_0(x)$. This implies $z \in\bigon_0(y_1)$ and $z \in\bigon_0(y_2)$. We assume w.l.o.g. that the oval containing the bigons connecting $x$ and $y_1$ and $y_2$ and $z$ comes before the oval containing the bigon connecting the other pairs. We then have $\mathrm{pre}(y_1,z)\not\equiv \mathrm{pre}(x,y_2)\pmod{2}$ and $\mathrm{pre}(y_2,z)= \mathrm{pre}(x,y_1)$. Because for a generator $g$, $I(g,g)$ does not depend on which side of the ovals the intersection points of $g$ are, we have $I(x,x)=I(y_2,y_2)=I(y_1,y_1)$. The claim follows.}
\end{itemize}
\end{proof}
\end{lemma}

\begin{theorem}
 Let $K$ be a knot with rectangular diagram $D$. The homology 
 of $(C'_{\Long}(D),\partial '_{\Long}(D))$ is $\widehat{HL}(K)\otimes  V^{\otimes(n-1)}$ (see Theorem \ref{MOSTh} for the definition of $V$). 
\begin{proof}
(sketch) We begin by extending our sign assignment for $(C_{\Long}(D),\partial_{\Long}(D))$ to a sign assignment for a modification of the long oval complex called $(C_{\mathrm{full}},\partial_{\mathrm{full}})$ that has the full $HL(K)\otimes  V^{\otimes(n-1)}$ (and not $\widehat{HL}(K))\otimes  V^{\otimes(n-1)}$) as homology. In \cite{moz} Section 2.3, the complex for $\widehat{HL}(K)\otimes  V^{\otimes(n-1)}$ is "extracted" from the more complicated complex for $HL$ in the same way that the long oval complex can be extracted from $(C_{\mathrm{full}},\partial_{\mathrm{full}})$.

 We then prove that there is, up to quasi-isomorphism, only one sign assignment $(C'_{\mathrm{full}},\partial'_{\mathrm{full}})$ with ${\partial'_{\mathrm{full}}}^2=0$. This is analogous with \cite{moz} Section 4.1. The existence of a sign assignment for $(C_{\mathrm{full}},\partial_{\mathrm{full}})$ with homology $HL(K)\otimes  V^{\otimes(n-1)}$ is a consequence of the analytical theory of link Floer homology (see \cite{knots} for the general theory and \cite{beliakova} for the special case of the long oval complex). Since our (extended) sign assignment verifies the condition ${\partial'_{\mathrm{full}}}^2=0$ (this is proved in the same way as Lemma \ref{d2iszero}), it must be quasi-isomorphic to the one coming from the analytical theory. Quasi-isomorphic complexes have identical homologies. The theorem follows by restricting the sign assignment on $\partial_{\mathrm{full}}$ to a sign assignment on $\partial_{\mathrm{long}}$ and noticing that because the homology of $(C'_{\mathrm{full}},\partial'_{\mathrm{full}})$ is $HL(K)\otimes  V^{\otimes(n-1)}$, the homology of $(C'_{\Long}(D),\partial '_{\Long}(D))$ will be $\widehat{HL}(K)\otimes  V^{\otimes(n-1)}$. 
\end{proof}
\end{theorem}
The short oval complex over ${\mathbb{Z}}/2{\mathbb{Z}}$ can now be extended to a complex with $\mathbb{Z}$ coefficients. We get a "signed" short oval complex by applying to the new long oval complex a tiny modification of the inductive construction we used in the $\mathrm{modulo\,}{2}$ case. Instead of Definition \ref{edgeDel}, we take:

\begin{definition}
\label{edgeDelSigned}
Let $p_1$ and $p_2$ be the intersection points that disappear at time $t=k$, assume $\{ p_1\}$ has a bigger Maslov grading than $\{ p_2\}$. 
For $x$ a generator of $C_{k}(D)$, let $\eta (x)$ be either $0$, if $x$ does not contain $p_2$, or, if $x$ does contain $p_2$, a generator identical to $x$ except that it contains $p_1$ instead of $p_2$.
$$\partial_k=\pi \circ(\partial_{k-1}-\partial_{k-1}\circ \eta \circ \partial_{k-1})\circ \iota$$
Where $\pi$ is the natural projection from $C_{k-1}(D)$ to $C_k(D)$ that sends generators containing $p_1$ or $p_2$ to zero and $\iota$ is the natural injection from $C_k(D)$ to $C_{k-1}(D)$. 
\end{definition}

 And instead of lemma 2.1 of \cite{beliakova} we use:
\begin{lemma}
\label{elimStep}
The complexes $C_{k}(D)$ and $C_{k-1}(D)$ of Definition \ref{edgeDelSigned} are homotopy equivalent.
\begin{proof}
This lemma can be proved either by "adding" signs to the proof in  \cite{beliakova} or by noticing that the generators containing a corner of the bigon can be matched in the sense of discrete Morse theory and that $C_{k}(D)$ is $C_{k-1}(D)$ with the matching collapsed (see \cite{kozlov}).
\end{proof}
\end{lemma}

 Our new short oval complex is not a sign assignment over the one defined in the previous section since some coefficients of the boundary map can be of absolute value greater than one.

\section{The Alexander grading}

In this section, we give a proof of Formula \ref{alexDef} for the Alexander grading of the long oval complex. Let us first formulate a technical lemma.

\begin{lemma}
\label{MinorSum}
The minimal number $l$ such that a permutation $\sigma$ is a product of $l$ transpositions is written $\si(\sigma)$. Let $M$ be a matrix of size $n$ with all entries in its first row and first column equal to 1.   Then, (\ref{lemma}) holds.
\begin{equation}
\label{lemma}
\det(M)=\sum_{1\le i,j\le n} (-1)^{i+j}\det(M_{i,j})
\end{equation}

Where $M_{i,j}$ denotes the minor of $M$ obtained by removing the row number $j$ and the column number $i$.
\begin{proof}
Let us denote the entries of $M$ by $m_{i,j}$ for $1\le i,j \le n$. Once we apply 
\begin{equation}
\label{detFormula}
\det(A)=\sum_{\sigma \in {\Sigma}_n} a_{1,\sigma (1)}...a_{n,\sigma (n)}(-1)^{\si(\sigma)}
\end{equation}
 to both sides of (\ref{lemma}), we get an equality between polynomials.
\begin{equation}
\sum_{\sigma \in {\Sigma}_n} m_{1,\sigma (1)}...m_{n,\sigma (n)}(-1)^{\si(\sigma)}=
\label{e3}
\end{equation}
\begin{equation}
\sum_{1<i,j\le n} (-1)^{i+j}
\sum_{\sigma \in {\Sigma}_n, \sigma(i)=j} (-1)^{\si(\sigma)+\sharp\{k:0<(k-i)(\sigma(k)-j)\}}\prod_{1\leq k\leq n, k\neq i}m_{k,\sigma(k)}
\label{e4}
\end{equation}
 We will prove it by comparing coefficients of monomials on both sides, while setting $m_{1,j}=1$ and $m_{i,1}=1$. We classify the monomials appearing
in (\ref{e3}) in two types: $m_{1,\sigma (1)}...m_{n,\sigma (n)}(-1)^{\si(\sigma)}$ for $\sigma \in \Sigma_n$ with $\sigma(1)=1$ and $m_{1,\sigma (1)}...m_{n,\sigma (n)}(-1)^{\si(\sigma)}$ with  $\sigma(1)\neq 1$. The sum of the monomials of the first type are exactly the monomials appearing in (\ref{e4}) coming from $\det(M_{1,1})$. The monomials of the second type all appear in (\ref{e4}) exactly three times, two times with the same sign as in (\ref{e3}), one time with opposite sign: (we use $\in$ to denote that a monomial is a summand of a polynomial written in canonical form)
 \begin{equation}
m_{2,\sigma(2)}...m_{n,\sigma(n)}(-1)^{\si(\sigma)}\in \det(M_{1,\sigma(1)}),\det(M_{\sigma^{-1}(1),1}),-\det(M_{\sigma^{-1}(1),\sigma(1)}).
\end{equation}
All monomials of (\ref{e3}) are now accounted for in (\ref{e4}). But monomials of the type:
\begin{equation}
\pm m_{1,\sigma(1)}...m_{i-1,\sigma(i-1)}
m_{i+1,\sigma(i+1)}...m_{n,\sigma(n)}\in \det(M_{i,\sigma(i)})
\end{equation}
for some $i\neq 1$ with $\sigma(1)\neq 1$ and $\sigma(i)\neq 1$, remain in (\ref{e4}). Setting $\overline{\sigma}$ to be equal to $\sigma$ composed with the transposition exchanging $i$ and $\sigma^{-1}(1)$, we notice that

\begin{equation}
(-1)^{\si(\sigma)+\sharp\{k:0<(k-i)(\sigma(k)-\sigma(i))\}+i+\sigma(i)}\prod_{1\leq k\leq n, k\neq i}m_{k,\sigma(k)}
(\in \pm\det(M_{i,\sigma(i)})) 
\end{equation}
is equal to
\begin{equation}
(-1)^{1+\si(\overline{\sigma})+\sharp\{k:0<(k-\sigma^{-1}(1))(\overline{\sigma}(k)-\sigma(i))\}+\sigma^{-1}(1)+\sigma(i)}
\prod_{1\leq k\leq n, k\neq \sigma^{-1}(1)}m_{k,\overline{\sigma}(k)}
\end{equation}

$$ \nonumber ({}\in \pm\det(M_{\sigma^{-1}(1),\sigma(i)})).$$

Therefore, those monomials cancel each other too.
\end{proof}
\end{lemma}
\begin{theorem}
 The Alexander grading on generators of the long oval complex on a grid diagram $D$ is given by the same formula as the Alexander grading on generators of the MOS complex:
$$
\label{AGrading2}
A(x)=\sum_{p\in x} a(p)- \frac{1}{2} \Bigl(\sum_{o\in O} 
a(o)\Bigr) - 
\left(\frac{n-1}{2}\right).
$$
Here, n is the complexity of $D$, $a(p)$ the winding number of $D$ around $p$ and $O$ one of the two types of punctures.
 
\begin{proof}
It follows from the analytical theory of knot Floer homology, that the \textit{relative} Alexander grading on pairs of generators must be defined in the same way in the two complexes  (see \cite{oms} and \cite{beliakova}). We prove the theorem by showing that the Euler characteristics of the long oval complex and the MOS complex are equal:
$$\chi(C_{\Long})=\chi(C_{\mathrm{MOS}}).$$ 
Since the Alexander polynomial is never null, this fixes the absolute Alexander grading.

Let $M$ be an $n\times n$ matrix, with entries $m_{i,j}$ equal to $t$ power the winding number of the knot around the grid diagram intersection point of coordinates $(i,j)$. The Euler characteristic of the MOS complex can be calculated by taking the determinant $M$ and multiplying by a constant $f(D)$ depending only on $D$. 
\begin{figure}[h]
\mbox{\epsfysize=7cm \epsffile{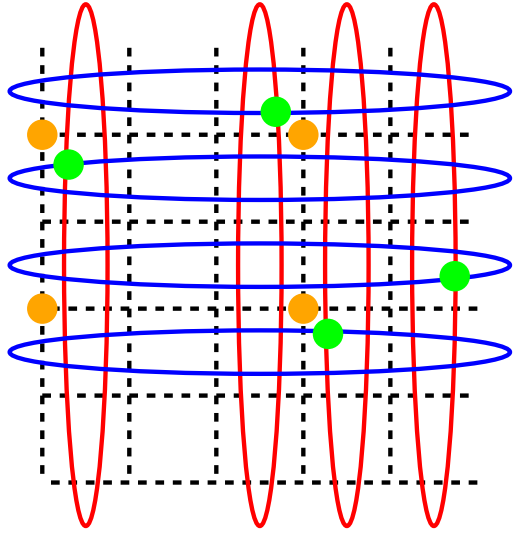}} 
\caption{{ A generator on a set of long ovals and the set of nearest intersection points on the grid.}}
\label{AlShift1}
\end{figure}

Drawing the diagram that gives rise to the long oval complex, as in Figure \ref{AlShift1}, we can associate sets of intersection points of the grid with oval complex generators. Working on the torus, we do this, as in Figure \ref{AlShift1}, by replacing intersection point of ovals by the nearest intersection point of the grid. It is easily checked (by pairing generators of opposite contributions, see Figure \ref{cancellation}) that the total contribution of the generators associated with sets having two points with the same vertical or horizontal coordinate is zero. But each set of intersection points of the grid with one point at most on each vertical or horizontal line of the grid is the associate of exactly one long oval generator.

\begin{figure}[h]
\mbox{\epsfysize=9cm \epsffile{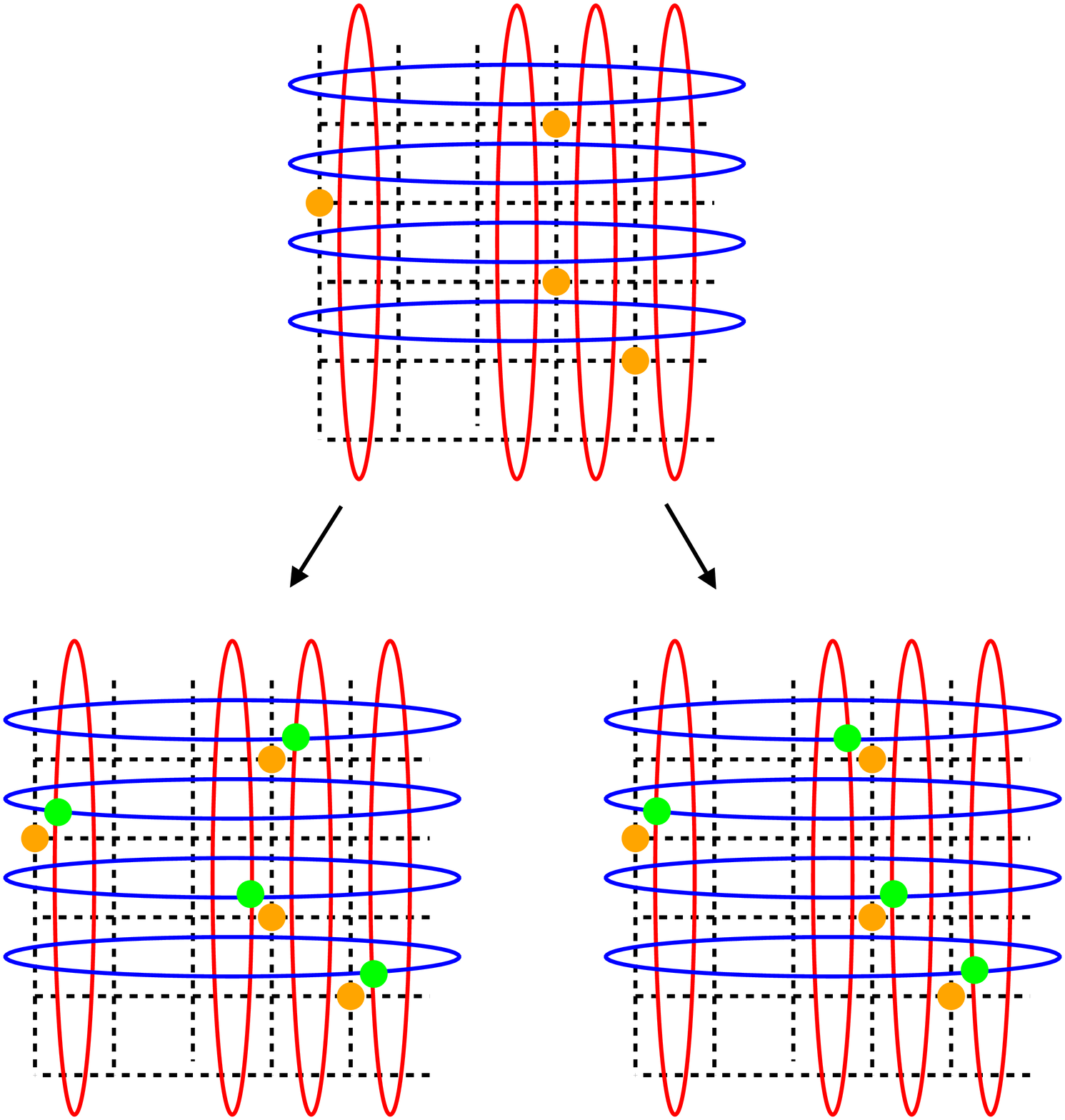}} 
\caption{{A set of $n-1$ intersection points with two points on the same column is associated here with two generators of the long oval complex. Since those generators have the same Alexander grading and since their Maslov gradings differ by one, they cancel each other in the Euler characteristic.}}
\label{cancellation}
\end{figure}
Let $S$ be any of those sets of $n-1$ points and $g$ its associated generator. The set $S$ misses exactly one row $j$ and one column $i$ of the grid diagram. The set $S$ can be seen as a monomial of the minor $M_{i,j}$, the point of coordinate $(k,l)$ in $S$ corresponding to the entry $m_{k,l}$ of $M$. The Alexander index of $g$ is exactly the value of the monomial of the minor $M_{i,j}$. Moreover, the sign of the monomial in $(-1)^{i+j}\cdot M_{i,j}$ is the Maslov index of $g$.
Collecting those facts, we get:
$$
\chi(C_{\Long})=\big(\sum_{m} {(-1)^m\cdot \sum_{a} {\mid C_{\Long}^{m,a}\mid \cdot t^a}}\big)\cdot f(D)\\
\nonumber=\sum_{1\le i,j\le n} (-1)^{i+j}\det(M_{i,j})\cdot f(D)
$$
Using Lemma \ref{lemma},
$$
=\det(M)\cdot f(D)=\chi(C_{\mathrm{MOS}}).
$$
\end{proof}
\end{theorem}


\section{The implementation}
For reasons of simplicity, we describe our program only for ${\mathbb{Z}}/2{\mathbb{Z}}$ coefficients. Passing to $\mathbb{Z}$ coefficient doesn't create any new problem.
Our program computes the knot Floer homology $\widehat{HL}(K)$ of a knot $K$ by constructing a short oval complex $(C_{\Short}(D(K)),\partial_{\Short}(D(K))$ and taking its homology. Our program therefore contains three main parts:
\begin{itemize}
\item{A function taking as input a knot (in braid representation usually) and giving as output a rectangular diagram with a set of ovals. The rectangular diagram and the ovals are chosen in order to minimize the number of generators of the short oval complex constructed from the ovals.} 
\item{A function listing the generators of the complex and calculating their Alexander and Maslov gradings.}
\item{A function that computes the boundary map. Since most of the running time of the program is spent by this function, it is heavily optimized.} 
\item{A function calculating the homology of the complex. This function gives in fact $\widehat{HL}(K)\otimes  V^{\otimes(n-1)}$. Extracting $\widehat{HL}(K)$ from it is easy.}
\end{itemize}
The first three parts are explained in the next three subsections.
The greatest task in the computation is to calculate the differential $\partial_{\Short}$. This differential, being a linear map, could be represented as a huge matrix $M_{\partial}$. However, because the differential preserves the Alexander grading and decreases the Maslov grading by one, the differential can be represented by a set of matrices of much more reasonable size, one for every pair of Alexander and Maslov gradings. 

Our short oval complex is a direct sum of complexes with uniform Alexander gradings.
$$
C_{\Short}=\bigoplus_{a} {C_{\Short}^a}
$$
Here, $C_{\Short}^a$ is the subcomplex of $C_{\Short}$ generated by generators of Alexander grading $a$. Since we compute an homology of the form $\widehat{HL}(K)\otimes  V^{\otimes(n-1)}$ with $V$ a vector space spanned by two vectors of different Alexander gradings, we don't even have to compute the whole complex $C_{\Short}$. For any set $A$ of $n-1$ Alexander gradings, the complexes $C_{\Short}^a$ can be deduced from the rest ($\bigoplus_{a\notin A} {C_{\Short}^a}$). For the Alexander gradings we choose to ignore, we don't have to compute the generators or the boundary map. Naturally, we choose to ignore the Alexander gradings that would be the most difficult to compute: those with many generators. Since the numbers of generators of the possible Alexander gradings is very unevenly distributed, a lot of time is saved.

\subsection{Generating a grid diagram}
\label{rectGen}
\begin{figure}[h]
\mbox{\epsfysize=120mm \epsffile{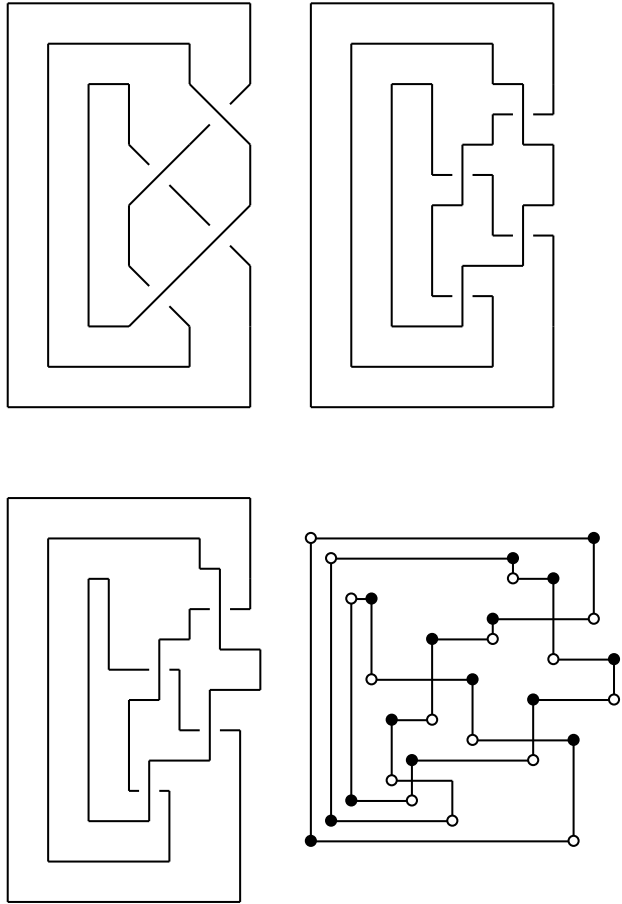}} 
\caption{{\ A representation by braid closure (top left) of the figure eight knot is transformed in four simple steps in a grid diagram (bottom right).}}
\label{rect}
\end{figure}
\begin{figure}[h]
\mbox{\epsfysize=120mm \epsffile{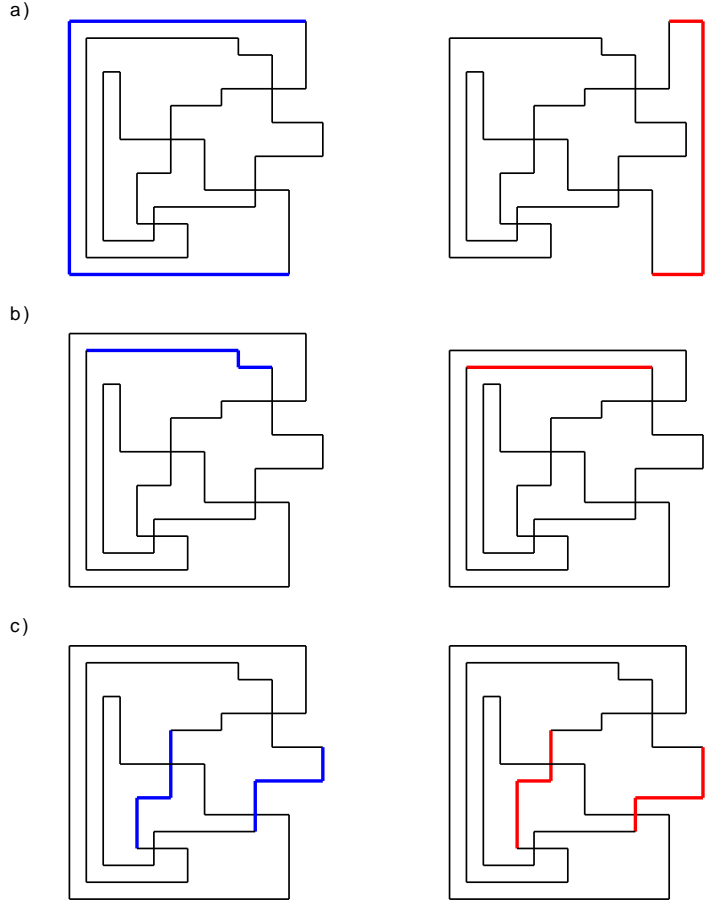}} 

\caption{{\ The three Cromwell moves:} (a) a cycling move, (b) a stabilization/destabilization move and (c) a castling move.}
\label{dyMoves}
\end{figure}

Generating a grid diagram representing a knot, either by hand or automatically, is not difficult (see Figure \ref{rect} for an example). However, the grid diagram obtained usually has much higher complexity than what is necessary. Therefore, our program simplifies knot diagrams before using them to construct complexes. This is done by exploring the set of equivalent rectangular diagrams during a fixed amount of time and picking the diagram of least complexity generated during this period. A set of simple moves (called Cromwell moves, see \cite{cromwell} or \cite{dy}) analogous to Reidemeister moves are used to this effect (see Figure \ref{dyMoves}). A cycling move is a circular permutation of the rows or column of a grid diagram. A stabilization move is the merging of two rows (or column) containing adjacent decorated squares accompanied by the deletion of the column (respectively the row) containing the two squares. A destabilization is simply the inverse of a stabilization. A castling move is an exchange of adjacent rows or columns of the grid diagram. (Castling moves are only possible if the decorations on the adjacent rows or columns are in a certain order.)
 The search is optimized in four ways:  
\begin{itemize}
\item{By caching (memorizing) already generated diagrams and exploring only diagrams of minimal complexity. This enables us to avoid exploring the same diagrams many times.}
\item{By considering as identical diagrams that can be obtained from one another by cycling moves. By reducing the space we have to explore, this makes caching much more beneficial. This is more or less equivalent to having the diagrams sitting on the torus instead of on the plane.}
\item{By exploring only diagrams of minimal complexity. In other words, moves increasing the complexity are not considered and each time a diagram of smaller complexity is discovered, the search begins anew from this diagram. This does not only gain time but also makes the memory consumption of the caching more reasonable. Moreover, searching in a monotonic way probably doesn't miss too many opportunities (see \cite{dy}).}
\item{By generalizing the stabilization move, so that simplifications can be discovered faster. For example, sequences of castlings in the same direction followed by a stabilization are considered as \emph{one} operation.}
\end{itemize}
Once we have a minimized grid diagram $D$ of complexity $n$, the next step is to embed this diagram in the plane and to choose $2n-2$ ovals on it. The generators of the complex $(C_{\Short}(D),\partial_{\Short}(D))$ that we want to calculate are sets of intersection point of ovals. In order to minimize the number of generators, we simply generate every possible embedding of the diagram and every choice of (shortened) ovals on them and pick the set of ovals that has the smallest number of intersection points. Of course, having a smaller number of intersection points of ovals does not guaranty having a smaller number of generators, but counting intersection points is easier than counting generators. 
\subsection{Listing generators}
\label{genGen}
Although constructing the generators is not very difficult, some care must be taken to compute their Alexander and Maslov gradings efficiently. To compute the Alexander gradings quickly, we tabulate $J(\{p\}-(\mathbb{O}+\mathbb{X})/2,\mathbb{X}-\mathbb{O})$
 beforehand for each intersection point $p$. Then for a given generator $g$, the Alexander grading is simply a sum of $n-1$ tabulated values and a precomputed constant.

The case of the Maslov grading is more difficult. The naive algorithm for computing $I(x,x)$ for a set $x\subset \mathbb{R}^2$, takes time $O(\vert x\vert^2)$. By using a divide and conquer approach (dividing $\mathbb{R}^2$ by horizontal and vertical lines that split $x$ in two for example), we can get an algorithm running in time $O(\vert x\vert \log {\vert x\vert})$. However, in this case, the constant is quite bad. Our solution consists in generating the generators in two stages. Each proto-generators\footnote{The proto-generators simply contains a pairing of the vertical ovals with the horizontal ones. They are therefore between $2^{n-1}$ and $4^{n-1}$ less numerous than the generators.} of the first stage giving raise to many generators in the second stage. This is done in a way that enables the worst part of the calculation, $I(x,x)$, to be done on the proto-generators and to be inherited by the generators. The rest of the calculation ($-I(x,\mathbb{O})-I(\mathbb{O},x)+I(\mathbb{O},\mathbb{O})$) uses the same tabulation method as the Alexander grading.

Since, the fiberedness and the genus of a knot only depend on the part of its homology with high Alexander grading, it is noteworthy that the set of generators $S$ of Alexander grading $>c$ can be listed in time $O(\vert S\vert \cdot n^3)$ (at least when $\vert S\vert>>n$). This comes from the fact that generators can be seen as matchings in a weighted bipartite graph between vertex representing the horizontal and vertical ovals. The weights can be chosen so that the Alexander grading of a generator is simply (up to a constant) the weight of the matching. Listing matchings with weights bigger than a constant is done by using a backtracking search that tests for cuts by the Hungarian algorithm for bipartite weighted matching (see \cite{kuhn}).

\subsection{The boundary map}
\label{bdmap}
Computing the boundary map means computing the entries $\partial_{x,y}$ of the matrix representing it as a function of the generators $x$ and $y$. When $\partial_{x,y}\neq 0$, we sometimes say that $x$ and $y$ are connected by the boundary map. We do not have a single method that always computes $\partial_{x,y}$ very quickly. However, $\partial_{x,y}$ is zero most of the time and, when $\partial_{x,y}$ is zero, we can usually determine that quickly. We therefore use two methods to compute $\partial_{x,y}$:
\begin{itemize}
\item{A fast method to check if there is a "domain" between two generators. The existence of a domain is a necessary condition for $\partial_{x,y}$ to be non-zero (see \ref{domains} below where the definition of domains is given).}
\item{A method that calculates $\partial_{x,y}$, which is inefficient, but which our program only calls when there is a domain between $x$ and $y$ (see \ref{chaining}).}

\end{itemize}
\subsubsection{Domains between generators}
\label{domains}
Let $D$ be a grid diagram on the plane with sets of vertical and horizontal ovals $\alpha$ and $\beta$. The ovals divide the plane in many bounded connected components and an unbounded component, which we call pieces. The bounded components are bigons, rectangles and triangles. \emph{Domains} are multisets (formal positive integral linear combinations) of pieces. For a domain $D$ and an intersection point $p$ of the ovals, the multiplicities in $D$ of the four pieces adjacent to $p$ listed in counter clockwise order starting from the upper right will be written $a_1(p),a_2(p),a_3(p),a_4(p)$. The corner index $c(p)$ of an intersection point $p$ is $a_1(p)+a_3(p)-a_2(p)-a_3(p)$. The sum of the pieces in an oval is called a \emph{periodic domain}.
\begin{lemma}
\label{periodic}
Given a set of transverse vertical and horizontal thin ovals in the plane, the only domains $D$ such that every intersection point of the curves has corner index zero are the sums of periodic domains and a multiple of the unbounded piece.
\begin{proof}
The result is obvious when the ovals have no intersection points. The lemma is then proved by induction on the number of intersection points.
\end{proof}
\end{lemma}
\begin{theorem}
\label{domTheo}
Given a collection of ovals on a grid diagram, the associated complex $(C,\partial)$ and two generators $x,y\in C$. If $\partial_{x,y}\neq 0$, there is a unique domain $D$ with the following properties:
\begin{itemize}
\item{The corner index at intersection points contained in $x$ but not in $y$ is 1.}
\item{The corner index at intersection points contained in $y$ but not in $x$ is -1.}
\item{All other intersection points have corner index zero.}
\item{Pieces of $D$ containing punctures have multiplicity zero in D.}
\end{itemize}
\begin{proof}
For long oval complexes, the existence and uniqueness of $D$ follow easily from the definition of the boundary map. For short oval complex, existence is proved by induction on the sequence of homotopies. 

Uniqueness of domain between generators in short oval complexes follows from the four conditions on $D$. Since the unbounded piece contains a puncture, the fourth condition implies that the multiplicity in $D$ of the unbounded piece is null. Therefore, using Lemma
\ref{periodic} and the four three conditions on $D$, we get uniqueness of $D$ up to addition of linear combination of periodic domains. However, at least one puncture is in only one periodic domain $p_0$ and all the periodic domains can be arranged in a sequence $p_0,\ldots ,p_{2n-2}$ such that $p_i$ and $p_{i+1}$ have one puncture in common. Using the fourth condition on $D$, this sequence enables us to show by induction on $i$ that the multiplicity of $p_i$ in the difference of two domains verifying the four conditions must be zero.
\end{proof} 
\end{theorem}
To determine the existence of a domain is therefore equivalent to solving in integers a system of linear equations and inequalities, the unknown being the multiplicities of the pieces in $D$. Luckily, the equations coming from the conditions on $D$ in Theorem \ref{domTheo} always have a unique solution, even over $\mathbb{Q}$. We can therefore decide efficiently if there is a domain between generators, by first solving a system of linear equations over $\mathbb{Q}$ and secondly checking that the solution obtained is non-negative and integral.

\subsubsection{Computing $\partial_{x,y}$}
\label{chaining}

The boundary map $\partial$ (for both long or short oval complexes) can be seen as a graph with vertices representing generators and edges representing the non-zero entries of the matrix $\partial$. Given an homotopy between long ovals and short ovals, the graph representing the short oval complex can be constructed by a sequence of modifications on the graph representing the long oval complex. The arguments used in this section are typical of algebraic Morse theory \cite{kozlov}.

Each time a pair of intersection points forming the corners of a bigon disappears (see Figure \ref{bigon2}), the graph is modified in the following way: 
\begin{itemize}
\item{All generators containing one of the disappearing intersection points are deleted.}
 \item{Let $a$ and $b$ be two generators connected by the bigon that disappears during the homotopy:
$$b\in \bigon_{0}(a).$$
 For every generators $x$ and $y$ with $\partial_{x,b}=1$ and $\partial_{a,y}=1$, we put one edge between $x$ and $y$ exactly if it wasn't there before.}   
\end{itemize}
This is a simple translation of Definition \ref{edgeDel} in graph theoretic language. Intuitively, it means that the edges of the graph of a short oval complex correspond to paths in the long oval complex graph. 

\begin{figure}[t]
\mbox{\epsfysize=140mm \epsffile{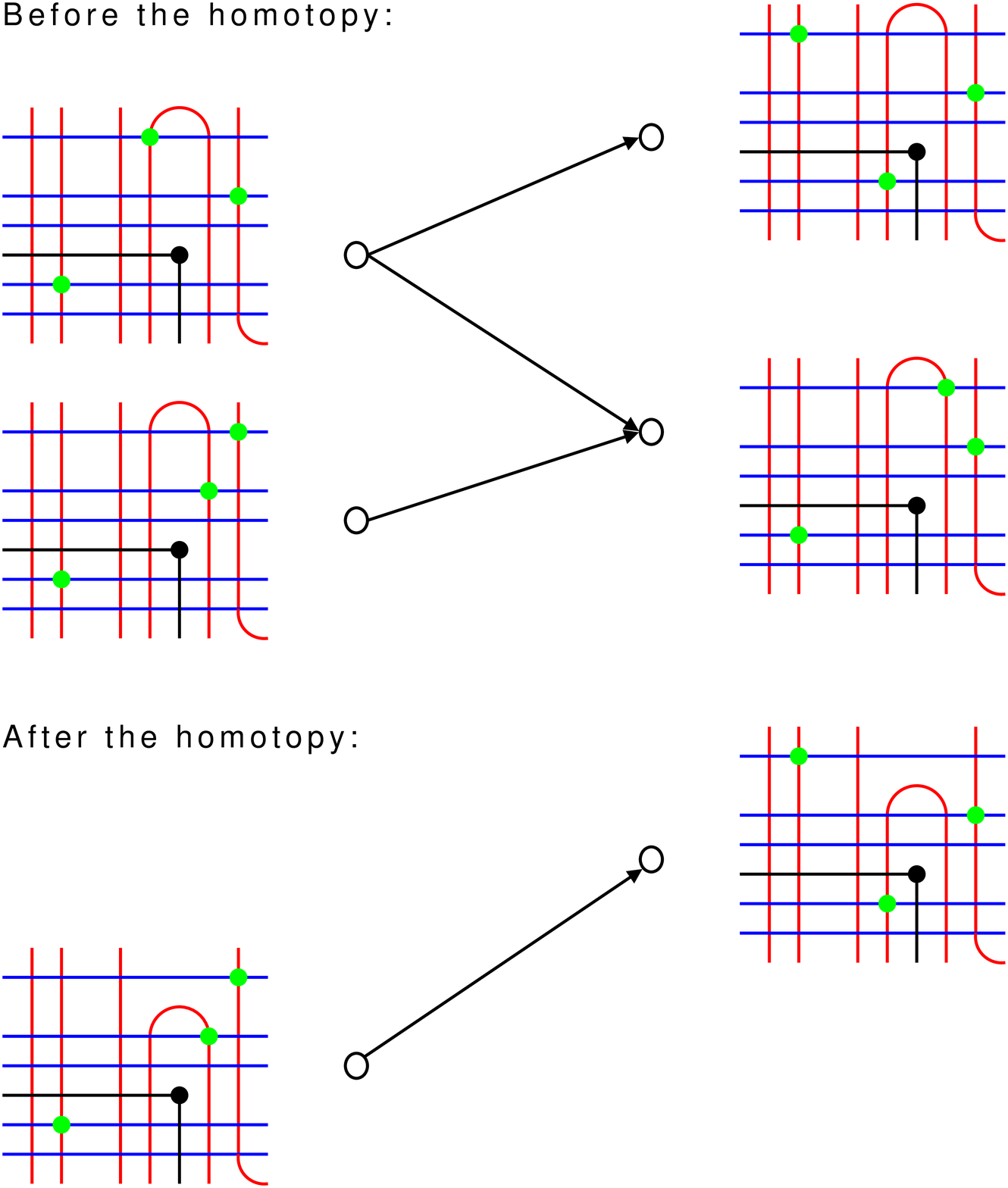}} 

\caption{{\ The effect of the disappearance of a bigon.  The small circles represent generators (a part of the corresponding generator is drawn beside each circle). The bipartite graph represents (a small part) of a boundary map between two consecutive Maslov gradings. The generators containing a disappearing intersection point disappear. We see that what was a path of length three between two generators, before the homotopy, becomes an edge between those two generators afterwards.}}
\label{path}
\end{figure}

Thinking in terms of graphs enables us to give a clearer condition for $\partial_{x,y}=1\pmod{2}$ in the short oval complex (see Figure \ref{path} for an illustration of the  effect of an homotopy on the graph representing the boundary map). 
\begin{theorem}
Let $(C_{\Long},\partial_{\Long})$ be a long oval complex. Let $(C_{\Short},\partial_{\Short})$ be a short oval complex obtained by an homotopy $(\alpha,\beta)_t$ for $t \in [0,l+1]$. Let $G$ be a graph representing $(C_{\Long},\partial_{\Long}).$ We associate to each edge of $G$, between generators connected by a bigon that disappears during the homotopy, the time at which it disappears. Then $x$ and $y$ are connected by the boundary map of the short oval complex exactly when the parity of the cardinality of the set of paths $W$ in $G$ verifying the following conditions is odd. Let $W$ be a sequence of edges $(e_1,...,e_k)$.
\begin{itemize}
\item{All the edges of $W$ are between generators of the same two Maslov gradings.}
\item{The path length $k$ is odd.}
\item{Each intermediate vertex on the path is adjacent to an edge of the path that disappears before it in the homology.}
\end{itemize}
\end{theorem}

Our program uses this combinatorial characterization of $\partial_{\Short}$ in a very straightforward way. It finds all possible paths in the long oval complex graph, by using a breadth first search, and it counts them. A simple caching of the generators explored is done to avoid exploring dead-ends many times. The point of this method is that the huge graph explored doesn't have to be stored at any time. It makes this method very efficient in terms of memory. Of course, the search is quite slow, but this method is not called often enough for speed to be a problem.

\subsection{Technical aspects of the programming}
Our program was written in Python\footnote{http://www.python.org/}. We used Idle and Eclipse \footnote{http://www.eclipse.org/} with PyDev \footnote{http://pydev.sourceforge.net/} as editing environment. Our program would have been much slower without the wonderful program called Psyco\footnote{http://psyco.sourceforge.net/} of Armin Rigo. Psyco is a kind of Just-In-Time compiler for Python which can apparently greatly accelerate (from two to a hundred time faster) almost any Python program. It enables programs in Python, a high level language, to run, even for the most computation intensive tasks, with speed comparable to programs in a low level language like C.

\bibliographystyle{custom}


\end{document}